\documentclass{article}

\usepackage{arxiv}
\usepackage{setspace}
\usepackage{lmodern}

\usepackage{textgreek}

\usepackage{fixltx2e}

\reversemarginpar
\usepackage{amsmath,amssymb,amsfonts,amsthm}
\usepackage{mathrsfs}% pour rajouter un format de lettres façon calligraphie en math mode.

\usepackage[locale = FR]{siunitx} % Pour gérer les unités
\DeclareSIUnit\vitesse{\meter\per\second}
\usepackage{eurosym}
\DeclareSIUnit{\octet}{o}

\usepackage{graphicx,array,tikz,multirow}
\usepackage{caption,subcaption} % permet de faire des subfigures (remplace le package subfig)
\usepackage{svg,float}
\usepackage{booktabs,paralist}
\newcolumntype{x}[1]{>{\centering\arraybackslash\hspace{0pt}}p{#1}}
\usepackage[section]{placeins}
\usepackage{hanging}

\usepackage{fancyhdr,emptypage} % garantit que les pages blanches avant les débuts de chapitres soient vraiment blanches (pas d'en-tête ni de pied de page)

\fancypagestyle{plain}{ %% Page chapitre, toc ...
    \fancyhead{}\fancyfoot[C]{\thepage}

}

\usepackage[Conny]{fncychap}

\usepackage[francais,nohints,tight]{minitoc}		% Mini table des matières, en français

\usepackage[nottoc]{tocbibind} % pour que la bibliographie apparaisse dans la table des matières (avec l'option pour que la table des matières elle-même n'apparaisse pas dans la table des matières).
\usepackage[titles]{tocloft}

\usepackage{textcomp} % rajoute des symboles 

\usepackage{titling}
\usepackage{lipsum} 
\usepackage{csquotes} % added 07/09/21

\usepackage{xspace}
\usepackage{afterpage}
\usepackage[textsize=footnotesize]{todonotes}

\usepackage{bookmark}
\usepackage{acronym}
\usepackage[nameinlink,french]{cleveref} % noabbrev
\Crefname{figure}{Fig.}{Figs.} % traduction des références aux figures/tables/équations
\crefname{figure}{fig.}{figs.}
\Crefname{equation}{Eq.}{Eqs.}
\crefname{equation}{eq.}{eqs.}
\Crefname{table}{Table.}{Tables.}
\crefname{table}{table.}{tables.}

\definecolor{color_ref}{rgb}{1.0, 0.13, 0.32} % couleur cite
\definecolor{color_link}{rgb}{0.0, 0.0, 1.0}
\definecolor{curcolor}{rgb}{0.0, 0.0, 1.0} % couleur des liens (bleu clair)
\definecolor{brightpink}{rgb}{1.0, 0.0, 0.5} %Bright pink
\definecolor{navyblue}{rgb}{0.0, 0.0, 1.0}
\hypersetup{
	colorlinks=true, % colore les liens au lieu de les encadrer
	pdfstartview=FitV, % ouvre le PDF de façon à ce qu'il prenne la taille verticale de l'écran
	urlcolor=brightpink, % choix de la couleur des liens URL
	linkcolor= navyblue, % choix de la couleur des liens internes (table des matières, etc.)
	citecolor=color_ref % choix de la couleur des liens de citations
}

\usepackage[utf8]{inputenc} % allow utf-8 input
\usepackage[T1]{fontenc}    % use 8-bit T1 fonts
\usepackage{hyperref}       % hyperlinks
\usepackage{url}            % simple URL typesetting
\usepackage{booktabs}       % professional-quality tables
\usepackage{amsfonts}       % blackboard math symbols
\usepackage{nicefrac}       % compact symbols for 1/2, etc.
\usepackage{microtype}      % microtypography
\usepackage{lipsum}
\usepackage{graphicx}
\usepackage{xcolor}
\graphicspath{ {./images/} }
\usepackage{amsmath}
\theoremstyle{definition}
\newtheorem{theorem}{Theorem}[section]

\newtheorem{definition}[theorem]{Definition}

\newtheorem{notation}[theorem]{Notations}

\newtheorem{remark}[theorem]{Remark}

\title{An adaptive method to solve multilevel multiobjective  linear programming problems}
\author{
 Mustapha Kaci$^{{\color{white}.}\bf 1}$ and Sonia Radjef$^{{\color{white}.}\bf 2}$\\
 $^{\bf 1,2{\color{white}.}}$Department of Mathematics, Signal Image Parole (SIMPA) Laboratory\\ University of Oran Mohamed Boudiaf USTO-MB, Oran, Algeria.\\
  {\color{blue}\texttt{kaci.mustapha.95@gmail.com} }$^{{\color{white}.}\bf 1}$\\
	{\color{blue}\texttt{soniaradjef@yahoo.fr} }$^{{\color{white}.}\bf 2}$
   }

\begin{document}
\fontfamily{ptm}\selectfont
\maketitle
\begin{abstract}
This paper is a follow-up to a previous work where we defined and generated the set of all possible compromises of multilevel multiobjective linear programming problems (ML-MOLPP). 
In this paper, we introduce a new algorithm to solve ML-MOLPP in which the adaptive method of linear programming is nested. First, we start by generating the set of all possible compromises (set of all non-dominated solutions). After that, an algorithm based on the adaptive method of linear programming is developed to select the best compromise among all the possible compromises achieved. This method will allow us to transform the initial multilevel problem into an ML-MOLPP with bounded variables. Then, apply the adaptive method which is the most efficient to solve all the multiobjective linear programming problems involved in the resolution process instead of the simplex method (It should be noted that the adaptive method is more efficient than the simplex method). Finally,  all the construction stages are carefully checked and illustrated with a numerical example.
\end{abstract}

\section{Introduction}

\par \hspace{0.95cm} In decision making, mathematical programming has long been restricted to problems having only one main objective, or a variety of objectives treated by striving to achieve them simultaneously. It is assumed that all the objectives are those of a single decision maker who controls all decision variables \cite{delh,prop3,ref555,ref2228,ref5,cm}. New decision making problems have emerged with a supplementary structure, these are problems in which several decision makers interact within a hierarchical structure, each seeking to optimize a multiobjective problem. Consequently, a new field of decision making theory appeared called multilevel mathematical programming \cite{prop6, prop5}. 
\par The principle of solving multilevel programming problems is that the first level decision maker (FLDM) sets its goal and/or decision, and then asks each subordinate level of the organization for its optimum, calculated in isolation. The lower-level decision-maker's decisions are then submitted and modified by the FLDM taking into account the overall benefits to the organization. The process continues until a satisfactory solution (also called satisfactory compromise) is found.
Multilevel programming problems have several applications in different branches such as supply chain management, network defense, planning, logistics, economics, government, autonomous institutions, agriculture, military, management, schools, hospitals, banks. Although most research on multilevel programming has focused on cases involving only two levels (called bilevel programming) \cite{baky,prop1,cite1,cc}, there are many programming problems that involve more than two levels \cite{baky2,kaci, prop2,refpuplishedpaper}.\\
Since the pioneering study of Bracken and McGill \cite{ref36, ref37}, several researchers have published monographs and literature surveys in which theoretical and methodological aspects of two-level optimization were discussed. Many approaches have been developed to solve the multilevel programming problems. For an extensive bibliography of these problems and their applications, see \cite{prop7, prop2,prop4,refpuplishedpaper, prop5}.
\par In 2003, a methodology using the e fuzzy goal programming (FGP) approach was proposed by Surabhi Sinha to solve the linear/non-linear multilevel decentralized programming problems (ML(D)PP) as described in \cite{prop8}. The method is subjective because it depends on the tolerance values given to the decision variables controlled by a decision-maker $(DM)$. It is costly over time due to the calculation of the bound of the decision variables and the increase in the number of constraints.
\par In their paper \cite{prop5}, Sinha and Sinha have introduced a linear programming approach applicable to any linear ML(D)PP. 
This is an improved version of S. Sinha’s method, as it does not depend on the subjectivity of the higher level $DM$, it provides a solution that is close to the ideal/optimal solution of each $DM$.
\par In his paper \cite{baky2}, I.A. baky extend the FGP approach introduced by  R.H. Mohamed \cite{prop10}  to solve ML-MOLPPs. The formulation of the FGP models begins by determining the fuzzy goals of the objectives by finding individual optimal solutions. They are then characterised by the associated membership functions; so also are the membership functions for vectors of fuzzy goals of the decision variables, controlled by $DM$s at the top levels. Moreover,  by introducing over- and under deviational variables and assigning aspiration level to each membership functions, they are transformed into fuzzy flexible membership goals. Then, the FGP is used to achieve the maximum value for each of the membership goals by minimizing their deviational variables and thus achieving the most satisfactory solution for all decision makers.
 \par In this study, we consider a ML-MOLPP where the objective  functions and the constraints are linear. We exploit the algorithm described by M. Kaci and S. Radjef in \cite{kaci} to generate the entire set of all the possible compromises. Then, we use the principle of the adaptive method to propose a new procedure for solving this last.
\par  The adaptive method of linear programming is a numerical constructive method described by R. Gabasov, F.M. Kirillova and O.I. Kostyukova in late 1980s \cite{Gabasov, Gabasov2}. It was generalized to develop many methods on piecewise linear programming, quadratic programming \cite{radjef}, optimal control,  multiobjective linear programming \cite{delh,ref5}, and multilevel multiobjective linear programming, this study.\\
The simplex method is a numerical method developed in 1947 for solving linear programming problems, where the feasible region $S$ is constructed from only linear constraints "Ax=b" and non-negativity constraints "$x\geq0$". If the linear program contains additional constraints of the form:
\begin{equation}\label{conn}
l\leq x\leq u,
\end{equation}
then, we say that we are dealing with linear programming problem with bounded variables. So, to solve this problem, we have to transform the  constraints ($\ref{conn}$) into linear constraints, which will increase the size of the matrix $A$.  Therefore, R. Gabasov develops the adaptive method to solve this kind of problem by handling the constraints ($\ref{conn}$) as such, and without modifying the matrix $A$. That made the adaptive method more efficient than the simplex method.\\
 The key idea of this study is to involve constraints of the form (\ref{conn}) in order to exploit the adaptive method. In other words, transform the problem (\ref{mlpp1}) into a multilevel multiobjective linear programming problem with bounded variables.
\par The adaptive method in all its versions is dedicated to the resolution of linear programming problems with bounded variables. So, the basic idea of the algorithm we’re building here is to determine the bounds of decision variables, then apply the latter instead of the simplex method. To achieve this, we first generate the set of all the possible compromises of the ML-MOLPP. As this one is not convex, we choose a convex subset called sorting set on which the resolution process will be performed. Finally, we define the components of the lower bound of decision variables as the minimum between the components of the corners of the sorting set, similarly, the upper bound vector will be defined using the maximum operator. 
\par In this paper, we propose a contribution in the area of multilevel multiobjective linear programming, it involves generalizing the adaptive method for solving  ML-MOLPP. First, we start by generating the set of all the non-dominated extreme points $\hat{N}^{dex}$ of all the multiobjective linear programming problems that constitutes the ML-MOLPP using the Yu and Zeleny's method \cite{YU}. Then, we generate the set of all  possible compromises $\hat{N}$ using the procedure described by M. Kaci and S. Radjef  in \cite{kaci}. After that, a sorting set will be chosen (a convex subset of $\hat{N}$) as a new feasible region and we put the bounds of all decision variables. Finally, we start the search for the satisfactory compromise from a chosen sorting set by solving $P-1$ standard linear multiobjective programming problems with bounded variables using the adaptive method algorithm for multiobjective linear programming problems described in section \ref{MA}, see also \cite{ref5}.
\par This paper is structured as follows:  First, the mathematical problem is formulated in the next Section. In Section \ref{section2}, we present some preliminary, that is to say, non-dominated solutions, and recall the adaptive method algorithm for solving multiobjective linear programming problems (MOLPP), as well as some necessary results and notations. In Section \ref{section3}, we build an algorithm that generalizes the adaptive method algorithm to solve an ML-MOLPP. After that, we illustrate the method with a numerical example  within the Section \ref{section4}. Finally, a conclusion is given in Section \ref{section5}.

\section{Problem formulation}\label{section1}
Consider a $P$-level multiobjective linear programming problem $(P\geq2)$, and  denote $DM_{p}$  the decision maker at $p^{\text{th}}$ level that has control over the decision variables $\overline{x}^{p}=x_{p1},\ldots, x_{pn_{p}}\in\mathbb{R}^{n_{p}}$, $p=1, \ldots, P$, where $x=(\overline{x}^{1}, \ldots,\overline{x}^{P})^{t}$, $n=n_{1}+\ldots+n_{P} $. \\Let  $k=k_{1}+\ldots+k_{P} $, we define
$$
\begin{array}{cccc}
   F_{p}:& \mathbb{R}^{n_{1}}\times\mathbb{R}^{n_{2}}\times \ldots \times\mathbb{R}^{n_{P}}&\longmapsto&\mathbb{R}^{k_{p}}\\
&x&\longmapsto&F_{p}(x)=c_{p}x
																																																														 \end{array},
$$
where

$$\begin{array}{cccccc}
c_{p}&=&\left(\begin{array}{c}c_{p1}\vspace{0.22cm}\\c_{p2}\\\vdots\\c_{pk_{p}}\end{array}\right)&=&
\left(\begin{array}{cccccccccc}c_{p1}^{11}&\ldots&c_{p1}^{1n_{1}}&c_{p1}^{21}&\ldots&c_{p1}^{2n_{2}}&\ldots&c_{p1}^{P1
}&\ldots&c_{p1}^{Pn_{P}}\vspace{0.22cm}\\
c_{p2}^{11}&\ldots&c_{p2}^{1n_{1}}&c_{p2}^{21}&\ldots&c_{p2}^{2n_{2}}&\ldots&c_{p2}^{P1
}&\ldots&c_{p2}^{Pn_{P}}\\
\vdots&&\vdots&\vdots&&\vdots&&\vdots&&\vdots\\
c_{pk_{p}}^{11}&\ldots&c_{pk_{p}}^{1n_{1}}&c_{pk_{p}}^{21}&\ldots&c_{pk_{p}}^{2n_{2}}&\ldots&c_{pk_{p}}^{P1
}&\ldots&c_{pk_{p}}^{Pn_{P}}\end{array}\right),
&p=1,\ldots,P
\end{array}
$$
and

$$\begin{array}{lll}
c_{pq}x&=&
c_{pq}^{1j}\overline{x}^{1}+c_{pq}^{2j}\overline{x}^{2}+\ldots+c_{pq}^{Pj}\overline{x}^{P},\hspace{0.4cm} p=1,\ldots,P,\hspace{0.3cm} q=1,\ldots,k_{p},\hspace{0.3cm} j=1,\ldots,n_{p}\vspace{0.2cm}\\
                  &=&c_{pq}^{11}x_{11}+\ldots+c_{pq}^{1n_{1}}x_{1n_{1}}+c_{pq}^{21}x_{21}+\ldots+c_{pq}^{2n_{2}}x_{2n_{2}}+\ldots+c_{pq}^{P1}x_{P1}+\ldots+c_{pq}^{Pn_{P}}x_{Pn_{P}}
\end{array}
$$
and
$$
\begin{array}{ccc}
c_{pq}^{pj}=c_{pq}^{p1},c_{pq}^{p2},\ldots,c_{pq}^{pn_{p}},& p=1,\ldots,P,& q=1,\ldots,k_{p}.
\end{array}
$$
The formulation of a $P$-level multiobjective linear programming problem is given as follows:

\begin{equation}\label{mlpp1}
\begin{array}{ll}
 \textbf{Level 1} & \\
& \underset{\overline{x}^{1}}{\max}\hspace{0.1cm} F_{1}(x)=\begin{array}{ll}\underset{\overline{x}^{1}}{\max}&\left(\begin{array}{c}c_{11}x\vspace{0.22cm}\\c_{12}x\\\vdots\\c_{1k_{1}}x\end{array}\right)\end{array},\vspace{0.25cm} \\
& \text{such that } \overline{x}^{2},\ldots,\overline{x}^{P} \text{solves }\\
 \textbf{Level 2} & \\
& \underset{\overline{x}^{2}}{\max}\hspace{0.1cm} F_{2}(x)=\begin{array}{ll}\underset{\overline{x}^{2}}{\max}&\left(\begin{array}{c}c_{21}x\vspace{0.22cm}\\c_{22}x\\\vdots\\c_{2k_{2}}x\end{array}\right)\end{array},\\
& \vdots\\
& \text{such that } \overline{x}^{P} \text{solves }\\
  \textbf{Level P}  & \\
& \underset{\overline{x}^{P}}{\max}\hspace{0.1cm} F_{P}(x)=\begin{array}{ll}\underset{\overline{x}^{P}}{\max}&\left(\begin{array}{c}c_{P1}x\vspace{0.22cm}\\c_{P2}x\\\vdots\\c_{Pk_{P}}x\end{array}\right)\end{array},
\end{array}
\end{equation}

subject to
$$
x\in S=\left\{x\in\mathbb{R}^{n} \hspace{0.15cm}:\hspace{0.15cm} Ax\leq b,\hspace{0.15cm} x\geq 0,\hspace{0.15cm} b\in\mathbb{R}^{m}\right\}.
$$
Where $S\neq\emptyset$ is the multilevel convex constraints feasible choice set, $m$ is the number of constraints, $k_{p}$ is the number of objective functions of the $DM_{p}$, $c_{pq}^{ij}$  are constants,  $A$ is $(m\times n)$-matrix and $b$ is a $m-$vector. We assume that $m<n$, and that the feasible region $S$ is a compact polyhedron  (a closed and bounded subset of $\mathbb{R}^{n}$).

%**********************************************************************************************************************************************************************SECTION 2******************************************************************

\section{Preliminaries}\label{section2}
\par In this section, we recall some important results on multilevel multiobjective linear programming and present the adaptive method which is a variant of the direct support method for solving multiobjective linear programming problems with bounded variables, we will also set some notations. For more details, see \cite{kaci,YU}.
\begin{notation}
For any $x\in S$, we use $x_{i}$ to indicate its $i^{\text{th}}$ coordinate, also called its $i^{\text{th}}$ component, and we define the following sets of indexes:
$$
\begin{array}{llll}
I=\left\{1,2,\ldots,m\right\},&J=\left\{1,2,\ldots,n\right\},&M=\left\{1,2,\ldots,n+m\right\},&J=J_{B}\cup J_{N},\\
J_{B}\cap J_{N}=\emptyset,  &\left|J_{B}\right|=m,&\left|J_{N}\right|=n-m.&
\end{array}
$$
This will allow us to write the partition of $x$ and $A$ as follows:
$$
\begin{array}{ll}
x=x(J)=(x_{j},j\in J)\hspace{0.25cm},& x=\left(\begin{array}{c}x_{N}\\x_{B}\end{array}\right),\vspace{0.2cm}\\
x_{N}=x(J_{N})=(x_{j}, j\in J_{N}),\hspace{0.25cm}&x_{B}=x(J_{B})=(x_{j},j\in J_{B})
\end{array}
$$
and
$$
A=(A_{N},A_{B}),  \, \, \,A_{N}=A(I,J_{N}),  \, \, \, A_{B}=A(I,J_{B}).
$$
The decision variable vectors  $x\in S$  will be considered as column vectors, $x^{t}$ refers to the transposition of $x$, which is a row vector. The $x>0$ and $x\geq0$ ratings indicate that all $x$ components are positive and not negative, respectively. For two vectors $x^{1}$ and $x^{2}$, the $x^{1}>x^{2}$ notation means $x^{1}-x^{2}>0$. The $x^{1}\geq x^{2}$, $x^{1}<x^{2}$, $x^{1}\leq x^{2}\ldots$etc notations should be interpreted accordingly. 

\end{notation}
\begin{remark}$\;$\\
\begin{itemize}
\item Let $p=1,\ldots, P$ and consider the $p^{\text{th}}-$problem which constitutes the multilevel problem \ref{mlpp1} as follows 
\begin{equation}\label{Pblm12022}
 \underset{\overline{x}^{p}}{\max} \hspace{0.2cm} F_{p}(x)=c_{p}x.
\end{equation}
Then, obviously the problem \ref{Pblm12022}  is not a real MOLPP. Indeed, the maximization occurs with respect to certain decision variables $\overline{x}^{p}=x_{p1},\ldots, x_{pn_{p}}$ only. 
\item
Since the direction of the objective function is not always "maximization", we could write "optimize". However, we choose the maximum operator to simplify the comprehension of this study.
\end{itemize}
\end{remark}
In the following, we will treat problems \ref{Pblm12022} as real multiobjective programming problems. So, consider the following multiobjective linear programming problem:
\begin{equation}\label{Pblm1}
     \underset{x\in S}{\max} \hspace{0.2cm} F_{p}(x)=c_{p}x. 
\end{equation}
\begin{definition}
 Let $Z_{p}$ be the $DM_{p}$'s criteria space  defined as follows:
$$
Z_{p}=\left\{z\in \mathbb{R}^{k_{p}}: z=c_{p}x,\hspace{0.15cm} x\in S\right\}.
$$
Then,
\begin{itemize}
\item For $z^{1},z^{2}\in Z_{p}$, we say that $z^{1}$ dominates $z^{2}$ if $z^{1}\geq z^{2}$.
\item For the feasible solutions $x^{1},x^{2}\in S$, we say that $x^{1}$ dominates $x^{2}$ if $c_{p}x^{1}\geq c_{p}x^{2}$.
\item A feasible solution $x\in S$ is non-dominated if it is not dominated by any other feasible point of $S$.
\end{itemize}
\end{definition}

\begin{notation}	$\;$\\
\begin{itemize}
\item $N_{p}$ denotes the set of all \textit{non-dominated} feasible solutions of  problem (\ref{Pblm1}).
\item  Since the feasible region $S$ is a compact polyhedron, then it has a finite number of vertices (also called extreme points). We denote the set of all extreme points of $S$ by $S^{dex}$.
\item We denote  the set of all non-dominated extreme points (non-dominated feasible solution that belongs to $S^{dex}$) of problem (\ref{Pblm1}) by: 
$$
N^{dex}_{p}=S^{dex}\cap N_{p}
$$
 and  $N^{dx}_{p}$ will denote an arbitrary element of $N^{dex}_{p}$.
\item
Let $\tilde{A}$ denote a $((n+m)\times n)-$matrix which the first $m$ rows correspond to the rows of the matrix $A$ and   the last rows correspond to the $n$ non-negativity constraints $\left(x_{ij}\geq 0,\hspace{0.15cm} i=1,\ldots,P , \hspace{0.15cm}j=1,\ldots, n_{i} \right)$,  $\tilde{b}$ be a $(n+m)-$vector with  $\tilde{b}_{i}=b_{i}$ for $i=1, \ldots, m$ and $\tilde{b}_{i}=0$ for $i=m+1,\ldots,m+n$, which means:
$$
\begin{array}{cc}
\tilde{A}=\left(\begin{array}{c}A\\-Id_{n}\end{array}\right),&\tilde{b}=\left(\begin{array}{c}b\\0\\0\\\vdots\\0\end{array}\right),
\end{array}
$$
 where $Id_{n}$ is the identity matrix of order $n$.
\end{itemize}
\end{notation} 
\begin{definition}

Let $Q\subseteq M$, $\tilde{A}_{Q}$ be the matrix derived from $\tilde{A}$ by deleting the rows which are not indexed in $Q$, similarly $\tilde{b}_{Q}$ is derived. Then:
\begin{itemize}
\item We call a facet of $S$, the set defined by:
$$
F(Q)=\left\{x\in S\hspace{0.15cm} :\hspace{0.15cm} \tilde{A}_{Q}x=\tilde{b}_{Q}\right\}.
$$
\item We call a non-dominated facet of $S$, the set defined by:
$$
N_{p}(Q)=N_{p}\cap F(Q).
$$
\end{itemize}
\end{definition}
\begin{definition}	$\;$\\	
\begin{itemize}
\item
 We define the set of all possible compromises (non-dominated feasible solutions) of a ML-MOLPP  as follows:
\begin{equation}
	\hat{N}=\bigcap_{p=1}^{P}{N_{p}}.
	\end{equation}
	\item
	Similarly, we define the set of all non-dominated extreme points of a ML-MOLPP as the set of common non-dominated extreme points  among all the sets $N_{p}$, that is:
	\begin{equation}\label{Ndex}
	\hat{N}^{dex}=\bigcap_{p=1}^{P}{N_{p}^{dex}}.
	\end{equation}
	\end{itemize}
\end{definition}

Assume that $ \hat{N}\neq\emptyset$, and define the following set of all subsets $Q$ of $M$ that corresponds to non-empty facets $F(Q)$ of the feasible region $S$ contained in $\hat{N}$:
	$$
	\hat{\mathcal{N}}=\left\{Q\subset M\hspace{0.22cm}:\hspace{0.2cm} F(Q)\neq\emptyset\hspace{0.22cm} \text{and}\hspace{0.22cm} F(Q)\subseteq\hat{N}\right\}.
	$$
	\begin{notation}
	The facet corresponding to a subset $Q$  of $\hat{\mathcal{N}}$ is called non-dominated facet for ML-MOLPP, and denoted by
	$$
	\hat{N}(Q)=\bigcap_{p=1}^{P}{N_{p}(Q)}.
	$$
	\end{notation}
	It is clear that for two subsets $Q_{1}$ and $Q_{2}$ in $\hat{\mathcal{N}}$ verifying $Q_{1}\subset Q_{2}$, we have $F(Q_{2})\subseteq F(Q_{1})$.  In order to eliminate all the facets $F(Q_{2})$ that are contained in bigger one $F(Q_{1})$, we define a new set $\tilde{\mathcal{N}}$ that contains the elements of $\hat{\mathcal{N}}$ which do not contain any other subset of $\hat{\mathcal{N}}$ as follows:
 
	$$
	\hat{\mathcal{N}}=\left\{Q\subset M\hspace{0.22cm}:\hspace{0.2cm} F(Q)\neq\emptyset\hspace{0.22cm} \text{and}\hspace{0.22cm} F(Q)\subseteq\hat{N}\right\}.
	$$

	\begin{definition}
A facet $F(Q)$, such that $Q\in\tilde{\mathcal{N}}$ is called \textit{sorting set} for ML-MOLPP and denoted by $\mathcal{SP}$.
	\end{definition}

\subsection{The adaptive method for solving multiobjective linear programming problems with bounded variables}  \label{MA}
 Consider the following multiobjective linear programming problem with bounded variables:
\begin{equation}\label{Pblm1221}
\left\{
\begin{array}{l}
     \underset{x}{\max} \hspace{0.2cm} c_{p}x \\ 	
      Ax\leq b\\
			l^{(p)}\leq x\leq u^{(p)}
\end{array}.
\right.
\end{equation}
After adding the slacks variables to the linear constraints $(Ax\geq b)$, we get $Bx=b$, where $B=\left(A,Id_{m}\right)$, and $Id_{m}$ is the identity matrix of size $m\times m$. To simplify the presentation, we will not change the notation of the decision variable vector $x$ and the matrix $c_{p}$. The use of the matrix $B$ means that $x, l^{(p)}, u^{(p)}$ are  $(n+m)-$vector  and $c_{p}=\left(c_{p},0_{k_{p}\times m}\right)$, where $0_{k_{p}\times m}$ is a matrix of size $(k_{p}\times m)$ with components equal to zero. Then, the problem (\ref{Pblm1221}) becomes:
\begin{equation}\label{Pblm122}
\left\{
\begin{array}{l}
     \underset{x}{\max} \hspace{0.2cm} c_{p}x \\ 	
      Bx=b\\
			l^{(p)}\leq x\leq u^{(p)}
\end{array},
\right.
\end{equation}
where $l^{(p)},u^{(p)}$ are a $(n+m)-$vectors such that for all $i\geq n+1$, $l^{(p)}_{i}=0$ and $u^{(p)}_{i}=\alpha$, $\alpha$ is a  larger positive number chosen in advance, and $|l^{(p)}|<\infty$, $|u^{(p)}|<\infty$. 
\par The adaptive method algorithm to solve the problem (\ref{Pblm122}) is described in the following scheme:

%%%%%%%%%%%%%%%%%%%%%%%%%%%%%%%%%%%%%%%%%%%%%%%%%%%%%%%%%%%%%%%%%%%%%%%%%%%%%%%%%%%%%%%%%%%%%%%%%

%%%%%%%%%%%%%%%%%%%%%%%%%%%%%%%%%%%%%%%%%%%%%%%%%%%%%%%%%%%%%%%%%%%%%%%%%%%%%%%%%%%%%%%%%%%%%%%%%%%%%
\begin{description}\label{moppalgo}
	\item[Step 1] Get a feasible solution $x^{0}$ of the problem (\ref{Pblm122}). If  $x^{0}$ exists then go to \textit{Step 2}, otherwise, stop the problem is infeasible.
	\item[Step 2] Get the solution $(y^{0},r^{0},v^{0},w^{0})$ of the following auxiliary  mono-objective linear programming problem:
	\begin{equation}
	\left\{\begin{array}{l}
	\begin{array}{l}y^{t}b-r^{t}c_{p}x^{0}-v^{t}l^{(p)}+w^{t}u^{(p)} \rightarrow \min                         \end{array}\\
	\begin{array}{l}y^{t}B+r^{t}c_{p}-v^{t}+w^{t}=s^{t}c_{p}                                         \end{array}\\
	\begin{array}{l}y\in\mathbb{R}^{m},\hspace{0.2cm}r\geq0,\hspace{0.2cm}w\geq0,\hspace{0.2cm}v\geq0\end{array}
	\end{array},\right.
	\end{equation}
	where $s^{t}=(1, \ldots, 1)\in\mathbb{R}^{k_{p}}$.
	\item[Step 3] Get the vector $\lambda_{p}=\left(r^{0}+s\right)$, and go to \textit{Step 4}.
	\item[Step 4] Elicit following mono-objective linear programming problem:
	\begin{equation}\label{Pblm12}
\left\{
\begin{array}{l}
     \begin{array}{l} \underset{x}{\max} \hspace{0.2cm} \lambda_{p}^{t}c_{p}x \end{array}\\
		 \begin{array}{l}Bx=b                                  \end{array}\\
		 \begin{array}{l}	l^{(p)}\leq x\leq u^{(p)}                \end{array}
		\end{array}
\right.
\end{equation}
\item[Step 5] Solve the linear programming problem (\ref{Pblm12}) using the adaptive method algorithm for mono-objective linear problems described in \cite{Gabasov}, to get a solution (non-dominated extreme point) of the problem (\ref{Pblm122}).
	\end{description}
	For more details, see \cite{delh, Gabasov2, Gabasov, ref5}.
	
%*****************************************************SECTION 04*********************************************************	*************************************************************************************************************************

	\section{The adaptive method for ML-MOLPP resolution}\label{section3}
	In this section, we  describe an algorithm that generalizes the adaptive method of linear multiobjective programming problems described in section \ref{MA}, in order to solve a ML-MOLPP  by choosing a satisfactory solution from a \textit{sorting set} $\mathcal{SP}$ chosen previously. 
	
	\subsection{Algorithm construction}\label{subsubsubs}
	Here is the detailed construction of the resolution algorithm:
	\begin{description}
	\item[Phase 1]$\;$\\ 
	\begin{enumerate}
	
	\item Generate the set of all possible compromises $\hat{N}$ of a ML-MOLPP using the method described in \cite{kaci}.
	\item
	\begin{itemize}
	\item Choose a sorting set $\mathcal{SP}\subset\hat{N}$ and assume that $
	|\hat{N}^{dex}|\geq2$,\hspace{0.15cm}$\mathcal{SP}\neq\hat{N}^{dex}$ and $\mathcal{SP}\neq\emptyset$.
	\item
	Define the set of extremes points of the set $\mathcal{SP}$ as follows:
	\begin{equation}
	\begin{array}{cllll}
	\mathcal{SP}^{dex}&:=&\hat{N}^{dex}\cap\mathcal{SP}&=&\left\{N^{dx_{1}},\ldots,N^{dx_{s}}\right\},
	\end{array}
	\end{equation}
 where $s$ is the cardinal of $\mathcal{SP}^{dex}$.
\end{itemize}
	\end{enumerate}
\item[Phase 2]	$\;$\\ 
	\begin{enumerate}
	\item  Put the bounds of each of the $n$ decision variables and slack variables as follows:
 \begin{equation}\label{bound}
\begin{array}{lll}
 l_{ij}=\underset{t\in\left\{1,\ldots,s\right\}}{\min}\left\{N^{dx_{t}}_{ij}\right\}\vspace{0.15cm},&
 u_{ij}=\underset{t\in\left\{1,\ldots,s\right\}}{\max}\left\{N^{dx_{t}}_{ij}\right\},
&  i=1,\ldots,P, \hspace{0.12cm}j=1,\ldots,n_{i},\\
l_{(P+1)j}=0,&u_{(P+1)j}=\alpha,&j=n+1,\ldots,n+m,
\end{array}
\end{equation}
where $\alpha$ is an arbitrary larger positive number chosen previously.
\item Put $l^{(1)}_{ij}=l_{ij}$,\hspace{0.12cm} $u^{(1)}_{ij}=u_{ij}$ for all  $i=1,\ldots,P$,\hspace{0.12cm} $j=1,\ldots,n_{i}$.												                      
\item Put $n_{P+1}=m$ and consider the following $n+m$ bounds of the decision and slack variables:
\begin{equation}\label{bound2}
\begin{array}{lr}
		     l^{(1)}_{ij}\leq x_{ij} \leq u^{(1)}_{ij},& i=1,\ldots,P+1,\hspace{0.12cm} j=1,\ldots,n_{i}.					
\end{array}
\end{equation}
\item Put $p=1$.
\item Define the new feasible region as follows:
$$
\begin{array}{c}
\mathcal{S}_{p}=\left\{x\in\mathcal{SP}\hspace{0.15cm}:\hspace{0.15cm}l^{(p)}_{ij}\leq x_{ij} \leq u^{(p)}_{ij}\right\}.
\end{array}
$$
\item 
\begin{itemize}
\item Solve the following multiobjective linear  programming problem with bounded variables:
\begin{equation}\label{pbsp}
\underset{x\in\mathcal{S}_{p}}{\max}\hspace{0.15cm}c_{p}x,
\end{equation}
\item Denote by $\stackrel{c}{x}^{p}$ the solution of the linear programming problem (\ref{pbsp}).
\end{itemize}
\item The $DM_{p}$ chooses $2n_{p}$ positive parameters $\textbf{\textit{l}}_{pj}$ and $\textbf{\textit{r}}_{pj}$ for $j=1,\ldots,n_{p}$. Then, define the new bounds of decision variables that are under his control as follows:
\begin{equation}\label{newbounds}
l^{(p)}_{pj}=\stackrel{c}{x}_{pj}^{p}-\textbf{\textit{l}}_{pj} \, \, \, \text {and} \, \, \,   u^{(p)}_{pj}=\stackrel{c}{x}_{pj}^{p}+\textbf{\textit{r}}_{pj},
\end{equation}
    such that
\begin{equation}\label{propriete1}
\begin{array}{l}
\stackrel{c}{x}_{pj}^{p}-\textbf{\textit{l}}_{pj}\geq l^{(1)}_{pj}\, \, \, \text {and} \, \, \,   \stackrel{c}{x}_{pj}^{p}+\textbf{\textit{r}}_{pj}\leq u^{(1)}_{pj}.
\end{array}
\end{equation}
\item $p=p+1$.
\item If $p>P$, then stop with a satisfactory compromise $\stackrel{c}{x}^{p}$ of the multilevel linear programming problem (\ref{mlpp1}). Otherwise,  go to Step 5.
	\end{enumerate}
\end{description}

		\subsection{Full algorithm}\label{amamlpp}
This is a simplified version (second formulation) of the resolution algorithm described in the previous section \ref{subsubsubs}:
	\begin{description}
  \item[Phase 1] $\;$\\
	\begin{enumerate}
	\item
	Apply the procedure described by M. Kaci and S. Radjef in \cite{kaci} to generate the entire compromises set $\hat{N}$ of ML-MOLPP.\vspace{0.15cm}
	\item Choose a sorting set $\mathcal{SP}\subset\hat{N}$ and elicit the set $\mathcal{SP}^{dex}$.\vspace{0.15cm}
	\end{enumerate}
  \item[Phase 2] $\;$\\
	\begin{description}
   \item[Step 1] Set $p=1$\vspace{0.15cm}.
	\item[Step 2] Formulate the constraints (\ref{bound2})\vspace{0.15cm}.
	\item[Step 3] Formulate the model (\ref{pbsp})\vspace{0.15cm}.
	\item[Step 4] $\;$\\
	\begin{itemize}
\item Solve the model (\ref{pbsp}) using the adaptive method algorithm described in section \ref{moppalgo}. \vspace{0.15cm}
\item Let $\stackrel{c}{x}^{p}$ denote the solution of the linear programming problem (\ref{pbsp}).
\end{itemize}
	\item[Step 5] Set $p=p+1$\vspace{0.15cm}.
	\item[Step 6] If $p>P$, then stop with a satisfactory compromise $\stackrel{c}{x}^{p}$ of the multilevel linear programming problem  (\ref{mlpp1}). Otherwise,  go to \textit{Step 7}.\vspace{0.15cm}
	\item[Step 7] Choose $\textbf{\textit{l}}_{p1},\ldots,\textbf{\textit{l}}_{pn_{p}}$ and $\textbf{\textit{r}}_{p1},\ldots,\textbf{\textit{r}}_{pn_{p}}$ such that the proprieties (\ref{propriete1}) holds in Step 8.\vspace{0.15cm}
	\item[Step 8] Set the $2n_{p}$ new  constraints (\ref{newbounds}), then go to Step 3.
\end{description}
\end{description}
%\end{description}
\subsection{Resolution process of ML-MOLPP with bounded variables}
\par
Multilevel programming problems are characterized by the presence of a hierarchical structure, that is, the presence of several decision-makers seeking a compromise. They are ranked by priority, the first decision-maker $DM_{1}$, the second decision-maker $DM_{2}$,$\ldots$, $P^{\text{th}}$ decision-maker $DM_{P}$(the last).\\
The $DM_{1}$ solves his multiobjective problem, depending on the solution he gets, he will impose constraints that ensure that any solutions that may appear in the following items do not leave a region he has defined. Then, he gives these new constraints to the $DM_{2}$ who will do the same; solve his problem and choose a new region (subset of the region established by the previous decision maker). So, the process continues until the solution is found. Therefore, the resolution process is as follows:\\
 We start by transforming the problem (\ref{mlpp1}) into an ML-MOLPP with bounded variables by defining the bounds $l^{(0)}$ and $u^{(0)}$ of the decision variables as indicated by the formulas (\ref{bound}).\\

 Consider the following positive component of the vectors $l^{(0)}$ and $u^{(0)}$ respectively:
$$
\begin{array}{ccc}
l^{(0)}_{11},\ldots,l^{(0)}_{1n_{1}},l^{(0)}_{21},\ldots,l^{(0)}_{2n_{2}},\ldots,l^{(0)}_{P1},\ldots,l^{(0)}_{Pn_{P}}\hspace{0.2cm} &\text{and}\hspace{0.2cm}& u^{(0)}_{11},\ldots,u^{(0)}_{1n_{1}},u^{(0)}_{21},\ldots,u^{(0)}_{2n_{2}},\ldots,u^{(0)}_{P1},\ldots,u^{(0)}_{Pn_{P}}
\end{array}
  $$
	\begin{description}
	\item[Iteration 1]$\;$\\
	\begin{itemize}
	\item[$\diamond$]
\par The $1^{\text{st}}-$decision maker begin by optimizing his objective function (the first objective function $F_{1}$) subject to the feasible region $S$ and constraints (\ref{conn}). That means, $DM_{1}$ solve the following multiobjective problem:
\begin{equation}
	\left\{
\begin{array}{l}
     \underset{x}{\max} \hspace{0.2cm}F_{1}(x)=c_{1}x \\
          Ax\leq b\\
					l^{(0)}\leq x\leq u^{(0)}
\end{array}.
\right.
\end{equation}
\item[$\diamond$] He gets his satisfactory solution denoted by
$
\stackrel{c}{x}^{1}=\left(\stackrel{c}{x}^{1}_{11}, \stackrel{c}{x}^{1}_{12}\ldots,\stackrel{c}{x}^{1}_{1n_{1}}, \stackrel{c}{x}^{1}_{21}\ldots,\stackrel{c}{x}^{1}_{P1},\ldots,\stackrel{c}{x}^{1}_{Pn_{P}}\right).
$
Then, impose a new constraints on the decision variables that are under his control:
$$
\stackrel{c}{x}^{1}_{11},\stackrel{c}{x}^{1}_{12},\ldots,\stackrel{c}{x}^{1}_{1n_{1}},
$$
by  choosing $2n_{1}$ reals numbers $l^{(1)}_{11}$, $l^{(1)}_{12}$, $\ldots$, $l^{(1)}_{1n_{1}}$, $u^{(1)}_{11}$, $u^{(1)}_{12}$, $\ldots$, $u^{(1)}_{1n_{1}}$ which verify 
$$
\begin{array}{ccccccc}
l^{(0)}_{11}\hspace{0.18cm}&\leq&\hspace{0.18cm}l^{(1)}_{11}\hspace{0.18cm}&<&\hspace{0.18cm} u^{(1)}_{11}\hspace{0.18cm}&\leq&\hspace{0.18cm} u^{(0)}_{11},\vspace{0.18cm}\\
l^{(0)}_{12}\hspace{0.18cm}&\leq&\hspace{0.18cm}l^{(1)}_{12}\hspace{0.18cm}&<&\hspace{0.18cm} u^{(1)}_{12}\hspace{0.18cm}&\leq&\hspace{0.18cm} u^{(0)}_{12},\\
\vdots&\vdots&\vdots&\vdots&\vdots&\vdots&\vdots\\
l^{(0)}_{1n_{1}}\hspace{0.18cm}&\leq&\hspace{0.18cm}l^{(1)}_{1n_{1}}\hspace{0.18cm}&<&\hspace{0.18cm} u^{(1)}_{1n_{1}}\hspace{0.18cm}&\leq&\hspace{0.18cm} u^{(0)}_{1n_{1}}.
\end{array}
$$
\item[$\diamond$]  Define the new constraints as follows:
\begin{equation}\label{lklk1}
\begin{array}{ccccc}
l^{(1)}_{11}&\leq& x_{11}&\leq& u^{(1)}_{11},\vspace{0.18cm}\\
l^{(1)}_{12}&\leq& x_{12}&\leq& u^{(1)}_{12},\\
\vdots&\vdots&\vdots&\vdots&\vdots\\
l^{(1)}_{1n_{1}}&\leq& x_{1n_{1}}&\leq& u^{(1)}_{1n_{1}}.
\end{array}
\end{equation}
\end{itemize}
\item[Iteration 2]$\;$\\
\begin{itemize}
\item[$\diamond$] 
\par The $2^{\text{nd}}-$decision maker solves his multiobjective problem subject to the feasible region $S$ and  constraints (\ref{conn} and \ref{lklk1}). That is, solve the following multiobjective linear problem:
\begin{equation}
	\left\{
\begin{array}{c}
   \underset{x}{\max} \hspace{0.2cm} F_{2}(x)=c_{2}x, \\
          Ax\leq b,\\
					\begin{array}{ccccc}
l^{\textbf{(1)}}_{11}&\leq& x_{11}&\leq& u^{\textbf{(1)}}_{11},\vspace{0.25cm}\\
l^{\textbf{(1)}}_{12}&\leq& x_{12}&\leq& u^{\textbf{(1)}}_{12},\\
\vdots&\vdots&\vdots&\vdots&\vdots\\
l^{\textbf{(1)}}_{1n_{1}}&\leq& x_{1n_{1}}&\leq& u^{\textbf{(1)}}_{1n_{1}},\vspace{0.25cm}\\
l^{(0)}_{21}&\leq &x_{21}&\leq& u^{(0)}_{21},\\
\vdots&\vdots&\vdots&\vdots&\vdots\\
l^{(0)}_{2n_{2}}&\leq& x_{2n_{2}}&\leq &u^{(0)}_{2n_{2}},\\
\vdots&\vdots&\vdots&\vdots&\vdots\\
l^{(0)}_{Pn_{P}}&\leq& x_{Pn_{P}}&\leq& u^{(0)}_{Pn_{P}}.
\end{array}
\end{array}
\right.
\end{equation}
\item[$\diamond$] The $DM_{2}$ obtains his satisfactory solution denoted by 
$$
\stackrel{c}{x}^{2}=\left(\stackrel{c}{x}^{2}_{11}, \stackrel{c}{x}^{2}_{12}\ldots,\stackrel{c}{x}^{2}_{1n_{1}}, \stackrel{c}{x}^{2}_{21}\ldots,\stackrel{c}{x}^{2}_{P1},\ldots,\stackrel{c}{x}^{2}_{Pn_{P}}\right).
$$
Then, define a new constraints on the decision variables that are under his control:
$$
\stackrel{c}{x}^{2}_{21},\stackrel{c}{x}^{2}_{22},\ldots,\stackrel{c}{x}^{2}_{2n_{2}},
$$
by choosing $2n_{2}$ reals numbers $l^{(2)}_{21}$, $l^{(2)}_{22}$, $\ldots$, $l^{(2)}_{2n_{2}}$, $u^{(2)}_{21}$, $u^{(2)}_{22}$, $\ldots$, $u^{(2)}_{2n_{2}}$ which verify   
$$
\begin{array}{ccccccc}
l^{(0)}_{21}\hspace{0.18cm}&\leq&\hspace{0.18cm}l^{(2)}_{21}\hspace{0.18cm}&<&\hspace{0.18cm} u^{(2)}_{21}\hspace{0.18cm}&\leq&\hspace{0.18cm} u^{(0)}_{21},\vspace{0.18cm}\\
l^{(0)}_{22}\hspace{0.18cm}&\leq&\hspace{0.18cm}l^{(2)}_{22}\hspace{0.18cm}&<&\hspace{0.18cm} u^{(2)}_{22}\hspace{0.18cm}&\leq&\hspace{0.18cm} u^{(0)}_{22},\\
\vdots&\vdots&\vdots&\vdots&\vdots&\vdots&\vdots\\
l^{(0)}_{2n_{2}}\hspace{0.18cm}&\leq&\hspace{0.18cm}l^{(2)}_{2n_{2}}\hspace{0.18cm}&<&\hspace{0.18cm} u^{(2)}_{2n_{2}}\hspace{0.18cm}&\leq&\hspace{0.18cm} u^{(0)}_{2n_{2}}.
\end{array}
$$
\item[$\diamond$] Define the new constraints as follows:
\begin{equation}\label{lklk2}
\begin{array}{ccccc}
l^{(2)}_{21}&\leq& x_{21}&\leq& u^{(2)}_{21},\vspace{0.18cm}\\
l^{(2)}_{22}&\leq& x_{22}&\leq& u^{(2)}_{22},\\
\vdots&\vdots&\vdots&\vdots&\vdots\\
l^{(2)}_{2n_{2}}&\leq& x_{2n_{2}}&\leq& u^{(2)}_{2n_{2}}.
\end{array}
\end{equation}
\end{itemize}
\item[Iteration 3]$\;$\\
\begin{itemize}
\item[$\diamond$] 
\par The $3^{\text{rd}}-$decision maker solves his multiobjective problem subject to the feasible region $S$ and  constraints (\ref{conn}, \ref{lklk1} and \ref{lklk2}). That is, solve the following multiobjective linear problem:
\begin{equation}
	\left\{
\begin{array}{c}
   \underset{x}{\max} \hspace{0.2cm} F_{3}(x)=c_{3}x,\\
          Ax\leq b,\\
					\begin{array}{ccccc}
l^{\textbf{(1)}}_{11}&\leq& x_{11}&\leq& u^{\textbf{(1)}}_{11},\vspace{0.25cm}\\
l^{\textbf{(1)}}_{12}&\leq& x_{12}&\leq& u^{\textbf{(1)}}_{12},\\
\vdots&\vdots&\vdots&\vdots&\vdots\\
l^{\textbf{(1)}}_{1n_{1}}&\leq& x_{1n_{1}}&\leq& u^{\textbf{(1)}}_{1n_{1}},\vspace{0.25cm}\\
l^{\textbf{(2)}}_{21}&\leq &x_{21}&\leq& u^{\textbf{(2)}}_{21},\vspace{0.25cm}\\
l^{\textbf{(2)}}_{22}&\leq &x_{22}&\leq& u^{\textbf{(2)}}_{22},\\
\vdots&\vdots&\vdots&\vdots&\vdots\\
l^{\textbf{(2)}}_{2n_{2}}&\leq& x_{2n_{2}}&\leq &u^{\textbf{(2)}}_{2n_{2}},\vspace{0.25cm}\\
l^{(0)}_{31}&\leq &x_{31}&\leq& u^{(0)}_{31},\\
\vdots&\vdots&\vdots&\vdots&\vdots\\
l^{(0)}_{Pn_{P}}&\leq& x_{Pn_{P}}&\leq& u^{(0)}_{Pn_{P}}.
\end{array}
\end{array}
\right.
\end{equation}
\item[$\diamond$]  The $DM_{3}$ obtains  his satisfactory solution denoted by 
$$
\stackrel{c}{x}^{3}=\left(\stackrel{c}{x}^{3}_{11}, \stackrel{c}{x}^{3}_{12}\ldots,\stackrel{c}{x}^{3}_{1n_{1}}, \stackrel{c}{x}^{3}_{21}\ldots,\stackrel{c}{x}^{3}_{P1},\ldots,\stackrel{c}{x}^{3}_{Pn_{P}}\right).
$$
Then,  defines new constraints that will be considered by the next decision-makers $DM_{4}$, $DM_{5}$,$\ldots$, $DM_{P}$.
\end{itemize}
\item[Iteration 4 $\longrightarrow$ P]
This process of resolution continues until a satisfactory compromise is reached in $P^{\text{th}}$-level.
\end{description}

\section{Numerical example}\label{section4}
Consider the following two-level multiobjective linear programming problem:
\begin{equation}
\begin{array}{ll}
 \textbf{Level 1} & \\
& \underset{\overline{x}^{1}=x_{1}}{\max}\hspace{0.1cm}\left(f_{11}(x)=2x_{1}+2x_{2}, f_{12}(x)=-\frac{1}{2}x_{1}+\frac{7}{25}x_{2}, f_{13}(x)=-\frac{1}{5}x_{1}+\frac{1}{2}x_{2}\right), \\
& \text{such that } \,  x_{2}  \, \, \text{solves }\\
 \textbf{Level 2} & \\
& \underset{\overline{x}^{2}=x_{2}}{\max}\hspace{0.1cm}\left(f_{21}(x)=x_{1}+3x_{2}, f_{22}(x)=-2x_{1}-x_{2}, f_{23}(x)=x_{2}\right).
\end{array}
\label{exemple-mlpp1}
\end{equation}
Subject to
$$
x\in S=\left\{x\in\mathbb{R}^{2}\hspace{0.15cm}:\hspace{0.15cm}Ax\leq b,\hspace{0.15cm} x\geq 0,\hspace{0.15cm} b\in\mathbb{R}^{7}\right\},
$$
where
$$
\begin{array}{ccc}
A=\left(\begin{array}{rrrrrrrrrrr}
-2& 1\\
-1& 2\\
 0& 1\\
 1& 0\\
-1&-2\\
 3&-4\\
 1&-2\\
\end{array}\right),&\hspace{0.3cm}
x=\left(\begin{array}{r}x_{1}\\ x_{2}\end{array}\right),&\hspace{0.3cm}b=\left(\begin{array}{r}3\\ 9\\ 6\\6 \\-9\\7\\2
\end{array}\right).
\end{array}
$$
\begin{description}
  \item[Phase 1]$\;$\\ 
	\begin{enumerate}
	\item The set of all the compromises of the ML-MOLPP (\ref{exemple-mlpp1}) is equal to: 
	$$
\begin{array}{ccc}
	\hat{N}&=&\mathcal{H}\left(N^{dx_{1}},N^{dx_{2}}\right)\cup\mathcal{H}\left(N^{dx_{2}},N^{dx_{3}}\right).
\end{array}
	$$
	Where $\mathcal{H}\left(.,.\right)$ denote the convex hull of any two points, the extremes points are
	$$
	N^{dx_{1}}=(6,6),\hspace{0.25cm} N^{dx_{2}}=(3,6),\hspace{0.25cm} N^{dx_{3}}=(1,5).
	$$
	\item
	 Choose the  sorting set $\mathcal{SP}$  equal to:
	$$
\mathcal{SP}=\mathcal{H}\left(N^{dx_{2}},N^{dx_{3}}\right).
  $$
Then, we get
	$$
	\mathcal{SP}^{dex}=\left\{N^{dx_{2}}=(3,6),N^{dx_{3}}=(1,5)\right\}.
	$$
	\end{enumerate}
	For more information about the calculation, refer to the numerical example in \cite{kaci}.
	\item[Phase 2] We put $p=1$, then:
	\begin{description}
	\item[Iteration 1] $\;$\\
	\begin{enumerate}
	\item
	The constraints (\ref{bound2}) are given by:
	$$
	\begin{array}{lr}
	l^{(1)}=(1,5),&
	u^{(1)}=(3,6).
	\end{array}
	$$
	\item Consider the following multiobjective programming problem:
	\begin{equation}\label{pbnv1}
	\left\{
\begin{array}{l}
    \underset{x}{\max} \hspace{0.2cm} c_{1}x \\
          Ax\leq b\\
					l^{(1)}\leq x\leq u^{(1)}
\end{array},
\right.
\end{equation}
where
$$
c_{1}=\left(\begin{array}{rr}
                 2&2\vspace{0.15cm}\\
							 -\frac{1}{2}&\frac{7}{25}\vspace{0.15cm}\\	
							 -\frac{1}{5}&\frac{1}{2}
             \end{array} \right).
						$$
	\item Using the adaptive method algorithm described in section \ref{moppalgo}, solve the multiobjective linear programming problem with bounded variables (\ref{pbnv1}). We start the adaptive method with the initial feasible solution:
	$$
	x^{0}=\left(2,\frac{11}{2}\right).
	$$
	Then, we get the compromise $\stackrel{c}{x}^{1}=\left(2,\frac{11}{2}\right)$ of $1^{st}-$level decision maker.
	\item Choose $\textit{\textbf{l}}_{1}=\textit{\textbf{r}}_{1}=0.5$.
	\end{enumerate}
	\item[Iteration 2] $\;$\\
	\begin{enumerate}
	\item Put
	$$
	\begin{array}{lllllll}
	l^{(2)}&=&x^{0}-(\textit{\textbf{l}}_{1},\textit{\textbf{l}}_{1})&=&\left(2,\frac{11}{2}\right)-(0.5,0.5)&=&\left(\frac{3}{2},5\right),\vspace{0.26cm}\\
	u^{(2)}&=&x^{0}+(\textit{\textbf{r}}_{1},\textit{\textbf{r}}_{1})&=&\left(2,\frac{11}{2}\right)+(0.5,0.5)&=&\left(\frac{5}{2},6\right).
	\end{array}
	$$
	\item
	Consider the following multiobjective linear programming problem:
	\begin{equation}\label{pbnv2}
	\left\{
\begin{array}{l}
     \underset{x}{\max} \hspace{0.2cm} c_{2}x \\
          Ax\leq b\\
					l^{(2)}\leq x\leq u^{(2)}
\end{array},
\right.
\end{equation}
where
$$
			c_{2}=\left(\begin{array}{rr}
                 1&3\vspace{0.15cm}\\
							 -2&-1\vspace{0.15cm}\\	
							  0&1
             \end{array} \right).
$$ 	
\item Take the initial feasible solution
	$$
	\left(\frac{5}{2},\frac{23}{4}\right).
	$$
	Then, use the adaptive method to solve the multiobjective linear programming problem (\ref{pbnv2}) to get a satisfactory compromise of the two-level multiobjective linear programming problem  (\ref{exemple-mlpp1}) equal to:
$$
\stackrel{c}{x}^{2}=\left(\frac{5}{2},\frac{23}{4}\right).
$$
\end{enumerate}
	\end{description}	
		\end{description}

\section{Conclusion}\label{section5}
\par 
In this study, the adaptive method for solving multiobjective linear problems is generalized to solve ML-MOLPP, we presented  an algorithm that searches for a single satisfactory compromise on a convex subset $\mathcal{SP}$ (sorting set) of the set of all possible compromises $\hat{N}$ of our choice (because $\hat{N}$ may be made up of several sorting sets). A detailed presentation of the method was followed by a numerical example.
\par
We show that the adaptive method continues to prove its importance in linear programming with regard to the diversity of problems in which it is applicable, since it was never been used to solve multilevel programming problems before namely the ML-MOLPP which was the main focus of this study, where the presence of bounded variables in the initial problem is no longer a requirement, they are established in the construction of the method.
\par Our main concern in this study was to give the construction of our approach based on the adaptive method that can be applied on the set of all possible compromises.  We will try in the near future to examine its robustness and report on the numerical results.

\bibliographystyle{unsrt}  
\bibliography{references}

\end{document}